\documentclass[11pt]{article}
\title
{Maker Can Construct a Sparse Graph on a Small Board}
\author{Heidi Gebauer
\thanks{Institute of
Theoretical Computer Science, ETH Zurich, CH-8092 Switzerland. Email:
gebauerh@inf.ethz.ch. }
}

\date{}
\usepackage{amsmath, amsfonts}
\usepackage{amsthm}
\usepackage{subfigure}
\usepackage[dvips]{graphicx}

\begin{document}
\bibliographystyle{plain}
\maketitle
\newtheorem{theo}{Theorem} [section]
\newtheorem{defi}[theo]{Definition}
\newtheorem{lemm}[theo]{Lemma}
\newtheorem{obse}[theo]{Observation}
\newtheorem{prop}[theo]{Proposition}
\newtheorem{coro}[theo]{Corollary}
\newtheorem{rem}[theo]{Remark}
\newtheorem{opprob}[theo]{Open Problem}

\newcommand{\whp}{{\bf whp}}
\newcommand{\prob}{probability}
\newcommand{\rn}{random}
\newcommand{\rv}{random variable}
\newcommand{\hpg}{hypergraph}
\newcommand{\hpgs}{hypergraphs}
\newcommand{\subhpg}{subhypergraph}
\newcommand{\subhpgs}{subhypergraphs}
\newcommand{\bH}{{\bf H}}
\newcommand{\cH}{{\cal H}}
\newcommand{\cT}{{\cal T}}
\newcommand{\cF}{{\cal F}}
\newcommand{\cD}{{\cal D}}
\newcommand{\cC}{{\cal C}}
\newcommand{\cP}{{\cal P}}
\newcommand{\cPG}{{\cal P}_{G}}

\newcommand{\ideg}{\mathsf {ideg}}
\newcommand{\lv}{\mathsf {lv}}
\newcommand{\nga}{n_{\text{game}}}
\newcommand{\avdeg}{\overline{\deg}}
\newcommand{\ed}{e_{\text{double}}}

\newcommand{\degb}{\deg_{B}}
\newcommand{\degm}{\deg_{M}}

\newcommand{\avd}{\bar{D}}

\newcommand{\remain}{\mathsf {rem}}

\newcommand{\dist}{\mathsf {dist}}
\newcommand{\pr}{\mathsf {Pr}}

\begin{abstract}
We study Maker/Breaker games on the edges of \emph{sparse} graphs. Maker and Breaker take turns in claiming previously unclaimed edges of a given graph $H$. Maker aims to occupy a given target graph $G$ and Breaker tries to prevent Maker from achieving his goal. We define a function $f$ on the integers and show that for every $d$-regular graph $G$ on $n$ vertices there is a graph $H$ with at most $f(d)n$ edges such that Maker can occupy a copy of $G$ in the game on $H$.
\end{abstract}

\section{Introduction}

We consider positional games played on edge-sets of graphs. Let ${\cP} = {\cP}(N) \subseteq 2^{E(K_{N})}$ be a graph property of $N$-vertex graphs, and let $H$ be a graph on the vertex set $V(H) = V(K_{N})$. The game $(E(H), {\cP})$ is played by two players, called Maker and Breaker, who take turns in claiming one previously unclaimed edge of $H$, with Maker going first. Following the standard notation we call $H$ the \emph{base graph} or the \emph{board}. Maker aims to occupy a graph having property $\cP$ and Breaker tries to prevent Maker from achieving his goal: Breaker wins if, after all edges of $H$ were claimed, Maker's graph does not possess ${\cP}$. 
A \emph{round} denotes a pair consisting of a Maker's move and the consecutive Breaker's move.

Let $G$ be a fixed graph on $n$ vertices. We consider the game where Maker's goal is to occupy a copy of $G$.
Formally, let ${\cP}_{G}$  denote the property that a graph contains $G$ as a subgraph and let $H = (V(H), E(H))$ be the board. The \emph{$G$-game} denotes the game $(E(H), {\cP}_{G})$. Note that Maker has only a chance to win if $|V(H)| \geq n$.

We show that if $G$ has maximum degree at most $d$ then there is a board $H$ where $|E(H)|$ is linear in $n$ such that Maker has a strategy to win the $G$-game on $H$.

\begin{theo} \label{theo:Makerwinsonsmallboard}
Let $G$ be a graph with maximum degree $d$ on $n$ vertices. Then there is a constant $c = c(d)$ and a graph $H$ with $|E(H)| \leq cn$ such that Maker has a strategy to occupy a copy of $G$ in the game on $H$.
\end{theo}

Feldheim and Krivelevich \cite{FK} showed that there are constants $c, c'$ depending on $d$ such that the following holds: if the board $H$ is the complete graph $K_{N}$ on $N = cn$ vertices then Maker can occupy a copy of $G$ in at most $c'n$ rounds (actually, they proved this statement also for the more general class of \emph{$d$-degenarate} graphs).
In our proof we will adopt many of their ideas, constructions and structures.

\paragraph{Notation}
Throughout this paper we will assume that Breaker starts the game.
Otherwise Maker can start with an arbitrary
move, then follow his strategy. If his strategy calls for something he occupied before he takes an
arbitrary edge; no extra move is disadvantegous for him. Accordingly, we slightly abuse notation and 
let a \emph{round} denote a pair consisting of a Breaker's move and the consecutive Maker's move.

Let $U, U' \in V(H)$. With $E_{H}(U,U')$ we denote the set of edges between $U$ and $U'$ in $H$.
Let $v \in V(G)$. The \emph{neighborhood} $N_{G}(v)$ denotes the set of vertices which are adjacent to $v$ in $G$. 
Let $u,v \in V(G)$. The \emph{distance} $\dist_{G}(u,v)$ \emph{between $u$ and $v$} denotes the number of edges in a shortest path in $G$ connecting $u$ and $v$. When there is no danger of confusion we sometimes omit the index $G$.
Adapting the notation of Feldheim and Krivelevich, the board graph along with the sets of Maker's and Breaker's claimed vertices is called a \emph{Game Position}.

\section{Proof of Theorem \ref{theo:Makerwinsonsmallboard}}

\subsection{Some Auxiliary Facts}

Let $G$ be a fixed $d$-regular graph on $n$ vertices.
Let $r = ed^{8}$. We will define an appropriate labeling $l$ which assigns to every vertex $v \in V(G)$ an element of $\{1, \ldots, r\}$. The \emph{level} of a vertex $v$ denotes $l(v)$. We first show that we can find a labeling with some particular properties.
\begin{lemm} \label{lemm:DivisionLevels}
There is a labeling $l$ of the vertices such that for every $u,v \in V(G)$ with $l(u) = l(v)$ we have $\dist(u,v) \geq 3$.
\end{lemm}
\emph{Proof:}
Suppose that we assign to each vertex $v \in V(G)$ a level in $\{1, \ldots, r\}$ uniformly and independently at random. Note that for any two vertices $u, v \in V(G)$, $\pr[l(u) = l(v)] = \frac{1}{r}$.
We will apply the symmetric version of the famous Lov\'{a}sz Local Lemma.
\begin{theo} (Symmetric version of the  Lov\'{a}sz Local Lemma.) \label{theo:symverloclem}
Let $A_{1}, A_{2}, \ldots, A_{n}$ be events in an arbitrary probability space. Suppose that each event $A_{i}$ is mutually independent of a set of all the other events $A_{j}$ but at most $k$, and that $\pr[A_{i}] \leq p$ for all $1 \leq i \leq n$. If
\begin{displaymath}
ep(k+1) \leq 1
\end{displaymath}
then $\pr[\bigwedge_{i=1}^{n} \bar{A}_{i}] > 0$.
\end{theo}
For every vertex $v \in V(G)$ let $A_{v}$ denote the event that $l(v) = l(w)$ for some vertex $w$ with $\dist(v,w) \leq 2$. Since $G$ has maximum degree $d$ there are at most $d + d^{2}$ vertices at distance at most 2 from $v$. So the probability $p = \pr[A_{v}]$ is at most $\frac{d + d^{2}}{r}$. Let $\{u_{1}, \ldots, u_{s}\}$ denote the set of vertices $u_{i}$ with $d(v, u_{i}) \leq 2$. 
Note that the event $A_{v}$ is completely determined by the values $l(v), l(u_{1}), l(u_{2}), \ldots, l(u_{s})$. Hence $A_{v}$ is independent of $\{A_{w} : \dist(v,w) \geq 5\}$. Applying the Local Lemma with $k = d + d^{2} + d^{3} + d^{4}$ yields that
\newline
$ep(k+1) \leq \frac{e(d + d^{2})(1 + d + d^{2} + d^{3} + d^{4})}{r} \leq \frac{ed^{3}d^{5}}{r} = 1$,
\newline
which concludes our proof. \hfill $\qed$

For every vertex $v$ let $d_{\leq 2}(v)$ denote the number of vertices at distance at most 2 from $v$ in $G$.
Suppose that $l$ is the labeling from Lemma \ref{lemm:DivisionLevels}.
Let $u,v$ be two vertices with $l(u) < l(v)$. We say that \emph{$u$ blocks $v$} if (i) $(u,v) \in E(G)$, \emph{or} (ii) there is a vertex $w$ with $l(w) < l(u) < l(v)$ such that $(w,u), (w,v) \in E(G)$.
We construct the directed graph $D$ on the vertex set $V(G)$ such that we draw an arc from $v$ to $u$ if and only if $u$ blocks $v$.
A vertex $u$ is called a \emph{descendant} of $v$ if there is a directed path from $v$ to $u$ in $D$.
\begin{obse} \label{obse:lengthpath}
Let $v \in V(G)$.
For every arc $(v,u)$ we have $l(v) > l(u)$. Moreover, every vertex $v$ has out-degree at most $d_{\leq 2}(v) \leq d + d(d-1) = d^{2}$ in $D$. Hence there are at most 
$(d^{2}) + (d^{2})^{2} + (d^{2})^{3} + \ldots + (d^{2})^{r-1} \leq (d^{2})^{r}$ descendants of $v$ in $D$.
\end{obse}

\subsection{Candidates and Candidate Schemes}

We first need some more notation.
The \emph{lower-level neighborhood} $N^{-}_{G}(v)$ of a vertex $v$ denotes the set $\{u \in N_{G}(v): l(u) < l(v)\}$. Accordingly, the
\emph{upper-level neighborhood} $N^{+}_{G}(v)$ denotes the set $\{u \in N_{G}(v): l(u) > l(v)\}$.

\paragraph{Construction of the board $H$} 
For every vertex $v \in V(G)$ we let $P(v)$ denote the set of descendants of $v$ in $D$.
Let $s_{d} := d^{5}2^{d+4}$. For every vertex $v \in V(G)$ we provide a set $S_{v}$ of $d(s_{d})^{2}|P(v)| + s_{d}$ vertices in $H$. So $V(H) := \cup_{v \in V(G)} S_{v}$. Moreover, for every $(u,v) \in E(G)$ we add an edge between every $a \in S_{u}$ and $b \in S_{v}$ in $H$. 
In other words, $E(H) := \{(a,b) : a \in S_{u}, b \in S_{v} \text{ such that $(u,v) \in E(G)$}\}$. 
Note that by Observation  \ref{obse:lengthpath}, $|P(v)| \leq d^{2r}$ and thus 
\begin{displaymath}
|E(H)| \leq |E(G)|(d(s_{d})^{2}d^{2r} + s_{d})^{2} \leq \frac{dn}{2} (d^{11 + 2ed^{8}} 2^{2d + 8} +d^{5}2^{d+4})^{2} 
\end{displaymath}
Note that $|E(H)|$ is linear in $n$.
As in \cite{FK}, to distinguish between vertices of $G$ and vertices of $H$ we mark the vertices of $H$ with a star.

During the game we will define for every vertex $v \in V(G)$ a subset $B_{v} \subseteq S_{v}$ with $|B_{v}| = s_{d}$.
We adopt the concepts of a \emph{candidate vertex} and a \emph{candidate scheme} from \cite{FK}, and state modified versions of Definition 2.1 - 2.3 (in \cite{FK}). 

\begin{defi} (Vertex candidate with respect to a specific edge) \label{defi:vertexcandidateforedge}
Let $H^{\star}$ be a position in $(E(H), \cPG)$, let $(u,v) \in E(G)$ with $l(u) < l(v)$ and let $\{u_{1}, \ldots, u_{t-1}\} = \{w \in N^{+}_{G}(u): l(w) < l(v)\}$. A vertex $x^{\star} \in S_{v}$ is called a \emph{candidate with respect to the edge $(u,v)$}, if 
\begin{itemize}
\item[(i)] $B_{u}, B_{u_{1}}, B_{u_{2}}, \ldots, B_{u_{t-1}}$ are already determined, and
\item[(ii)] for every choice of vertices $b^{\star}_{1} \in B_{u_{1}}, b^{\star}_{2} \in B_{u_{2}}, \ldots, b^{\star}_{t-1} \in B_{u_{t-1}}$ we have
    \begin{displaymath}
   \frac{ |\{b^{\star} \in B_{u} : \text{Maker claimed $(b^{\star}, b^{\star}_{1}), (b^{\star}, b^{\star}_{2}), \ldots, (b^{\star}, b^{\star}_{t-1}), (b^{\star}, x^{\star})$ in $H^{\star}$}\}|}{|B_{u}|} \geq \frac{1}{t2^{t}}
    \end{displaymath}
\end{itemize}
\end{defi} 

\begin{defi} (Vertex candidate) \label{defi:vertexcandidate}
Let $H^{\star}$ be a position in $(E(H), \cPG)$, let $v \in V(G)$ and let $x^{\star} \in B_{v}$. We call $x^{\star}$ a \emph{candidate} if for every $u \in N^{-}_{G}(v)$ $x^{\star}$ is a candidate with respect to $(u,v)$.
\end{defi}

Let $v_{1}, v_{2}, \ldots, v_{n}$ be an ordering of the vertices in $V(G)$ where $l(v_{1}) \leq l(v_{2}) \leq \ldots, \leq l(v_{n})$.

\begin{defi} \label{defi:candidatescheme}
Let $H^{\star}$ be a position in $(E(H), \cPG)$ and suppose that $B_{v_{1}}, B_{v_{2}}, \ldots, B_{v_{n}}$ are all determined. We say that $(B_{v_{1}}, B_{v_{2}}, \ldots, B_{v_{n}})$ form a \emph{candidate scheme} if every $x^{\star} \in B_{v_{1}} \cup B_{v_{2}} \cup \ldots \cup B_{v_{n}}$ is a candidate.
\end{defi}

The next lemma is a slight adaptation of Lemma 2.1 in \cite{FK}.
\begin{lemm} (Feldheim, Krivelevich) \label{lemm:CandidateSchemeImpliesGCopy}
Let $H^{\star}$ be a position in $(E(H), \cPG)$ and let $(B_{v_{1}}, B_{v_{2}}, \ldots, B_{v_{n}})$ be a candidate scheme. Then Maker's graph contains a copy of $G$.
\end{lemm}
A proof of this lemma is given in \cite{FK}.

\subsection{Dividing the $G$-Game into Subgames}

We first need some more notation.
Let $v \in V(G)$ and let $x^{\star} \in S_{v}$. By a slight abuse of notation we define the \emph{level} $l(x^{\star})$ of $x^{\star}$ to be $l(v)$. 
Let $H^{\star}$ be a position in the game $(E(H), {\cP}_{G})$. 
We call a vertex $x^{\star} \in S_{v}$ \emph{touched} in $H^{\star}$ if Maker or Breaker claimed an edge of the form $(y^{\star}, x^{\star})$ with $l(y^{\star}) < l(x^{\star})$.
Accordingly, we call $x^{\star}$ \emph{untouched} in $H^{\star}$ if $x^{\star}$ is not touched in $H^{\star}$.
We say that a vertex is $v \in V(G)$ \emph{completed} in $H^{\star}$ if every vertex in $B_{v}$ is a candidate.
We say that $v$ is \emph{ready} if every vertex $u$ with $(v,u) \in E(D)$ is completed.
When there is no danger of confusion we will sometimes omit mentioning $H^{\star}$ explicitly.

At the beginning of the game we set $B_{v} := S_{v}$ for every vertex $v \in V(G)$ with out-degree zero in $D$. Hence $v$ is already completed (and thus also ready). 

We now consider some smaller games separately and then show how Maker can divide the $G$-game into a combination of those smaller games. 
\begin{defi} \label{defi:addnewvertextogame}
Let $H^{\star}$ be a position in $(E(H), {\cPG})$, let $v \in V(G)$ be ready in $H^{\star}$ and let $x^{\star} \in S_{v}$. 

 For every $u \in N^{-}(v)$ the game $G_{u, x^{\star}}$ is defined as follows. The board consists of the set of edges $E_{H}(B_{u}, \{x^{\star}\})$; Maker's goal is to achieve that $x^{\star}$ becomes a candidate with respect to the edge $(u,v)$.
\end{defi}

We now extend this game to a game where not only $x^{\star}$  but several vertices of $S_{v}$ should become a candidate.

\begin{defi} \label{defi:addseveralverticestogame}
Let $H^{\star}$ be a position in $(E(H), {\cPG})$, let $v \in V(G)$ be ready in $H^{\star}$ and suppose that $B_{v}$ has been determined. The game $G_{v}$ is defined as follows. The board is $\cup_{u \in N_{G}^{-}(v)} E_{H}(B_{u}, B_{v})$; Maker's goal is to win $G_{u, x^{\star}}$ for \emph{every} $u \in N^{-}_{G}(v)$, $x^{\star} \in B_{v}$.
\end{defi}

The next lemma is an adaptation of Lemma 2.2 in \cite{FK}.
\begin{lemm} \label{lemm:adaptschemeextensionforsinglevert}
Let $H^{\star}$ be a position in $(E(H), {\cPG})$, let $v \in V(G)$ be ready and let $x^{\star} \in S_{v}$ be untouched in $H^{\star}$. Then, for every $u \in N^{-}_{G}(v)$ Maker has a strategy to win the game $G_{u, x^{\star}}$.
\end{lemm}
Let $N^{-}_{G}(v) = \{u_{1}, \ldots, u_{t}\}$ and note that if $v$ is ready then every $u_{i}$ and every 
$w \in N^{+}(u_{i})$ with $l(w) < l(v)$ is completed.
Then Lemma \ref{lemm:adaptschemeextensionforsinglevert} follows directly from the proof of Lemma 2.2 in \cite{FK}. For completeness we restate the core of the proof.

\emph{Proof of Lemma \ref{lemm:adaptschemeextensionforsinglevert}:}
We need some notation. Let $F = (V(F), E(F))$ be a hypergraph, i.e., $E(F)$ is a subset of the power set $2^{V(F)}$. In a \emph{positional game on $F$} Maker and Breaker alternately claim an unclaimed vertex of $V(F)$ until all vertices are claimed. We will use the following result by Alon, Krivelevich, Spencer and Szab\'{o} \cite{AKSS}, extending a previous result by Sz\'{e}kely \cite{S}.
\begin{theo} (Alon, Krivelevich, Spencer, Szab\'{o}) \label{theo:discrepancygame}
Let $F$ be a hypergraph with $X$ hyperedges, whose smallest hyperedge contains at least $x$ vertices. In a positional game on $F$ Maker has a strategy to claim at least $\frac{x}{2} - \sqrt{\frac{x \ln(2X)}{2}}$ vertices of each hyperedge.
\end{theo}

To prove Lemma \ref{lemm:adaptschemeextensionforsinglevert} we fix a $v \in V(G)$ and an untouched $x^{\star} \in S_{v}$. If $N^{-}_{G}(v) = \emptyset$ then the claim is clearly true. Otherwise let $u \in N^{-}_{G}(v)$. Let $\{u_{1}, \ldots, u_{t-1}\} = \{w \in N^{+}_{G}(u): l(w) < l(v)\}$ and recall that $t \leq d$. Note that since $v$ is ready we have that every $u_{i}$ is completed and therefore $B_{u_{1}}, \ldots, B_{u_{t-1}}$ are all candidates.

If $t=1$ then Maker can win the game $G_{u, x^{\star}}$ by connecting $x^{\star}$ to half of the vertices of $B_{u}$.
Suppose that $t \geq 2$. We can express $G_{u, x^{\star}}$ as a positional game on a hypergraph $F$: Let $V(F)$ be the set of all edges in $E_{H}(B_{u}, x^{\star})$; and for every combination of $b^{\star}_{1} \in B_{u_{1}}, b^{\star}_{2} \in B_{u_{2}}, \ldots, b^{\star}_{t-1} \in B_{u_{t-1}}$ add to $E(F)$ the hyperedge $e_{b^{\star}_{1}, b^{\star}_{2}, \ldots, b^{\star}_{t-1}}$ consisting of all edges of the form $(b^{\star}, x^{\star})$ where $b^{\star} \in B_{u}$ is connected to $b^{\star}_{1}, b^{\star}_{2}, \ldots, b^{\star}_{t-1}$ in Maker's graph.

We have $|E(F)| \leq (s_{d})^{t-1}$.
\newline
Note that $|e_{b^{\star}_{1}, b^{\star}_{2}, \ldots, b^{\star}_{t-1}}| =
|\{ b^{\star} \in B_{u} : \text{Maker claimed $(b^{\star}, b^{\star}_{1}), (b^{\star}, b^{\star}_{2}), \ldots, (b^{\star}, b^{\star}_{t-1})$} \}|$. Since $B_{u_{1}}, \ldots, B_{u_{t-1}}$ are all candidates we have 
\begin{displaymath}
|e_{b^{\star}_{1}, b^{\star}_{2}, \ldots, b^{\star}_{t-1}}| \geq \frac{|B_{u}|}{2^{t-1}(t-1)} \geq \frac{s_{d}}{2^{t-1}(t-1)}
\end{displaymath}
By Theorem \ref{theo:discrepancygame} Maker has a strategy to claim at least 
\begin{equation} \label{eq:applicationofdiscrepancytheorem}
\frac{s_{d}}{2^{t}(t-1)} - \sqrt{\frac{s_{d}}{2^{t}(t-1)} \ln(2(s_{d})^{t-1})}
\end{equation}
vertices in every hyperedge.
A careful calculation (details can be found in \cite{FK}) yields that the expression in \eqref{eq:applicationofdiscrepancytheorem} is at least 
$\frac{s_{d}}{2^{t}t}$. \hfill $\qed$

Lemma \ref{lemm:adaptschemeextensionforsinglevert} allows to show the following corollary, a similar version of which has already been stated in \cite{FK}.
\begin{coro} \label{coro:winningstrategyforgameononevert}
Let $H^{\star}$ be a position in $(E(H), {\cPG})$ and let $v \in V(G)$ be ready in $H^{\star}$. Suppose that $B_{v}$ has been determined and that every vertex in $B_{v}$ is untouched in $H^{\star}$. Then Maker has a strategy to win the game $G_{v}$.
\end{coro} 
\emph{Proof:} Consider the following strategy for Maker. Suppose that Breaker claims an edge $(b^{\star}, x^{\star})$ with $b^{\star} \in B_{u}$ and $x^{\star} \in B_{v}$ for some $u \in N^{-}_{G}(v)$ (note that due to the construction of the board every claimed edge is of this form). Then Maker responds in the game $G_{u, x^{\star}}$ 

Since the boards of the games $G_{u, x^{\star}}$ are pairwise disjoint Maker can treat each game $G_{u, x^{\star}}$ separately. Thus Lemma \ref{lemm:adaptschemeextensionforsinglevert} yields a winning strategy for Maker in $G_{v}$. \hfill $\qed$

The next observation shows that the game $G_{v}$ is finished after not too many rounds.
\begin{obse} \label{obse:gamelengthsinglevertex}
Let $v \in V(G)$. The board size of the game $G_{v}$ is bounded by $|N^{-}_{G}(v)| (s_{d})^{2} \leq d s^{2}{d}$, hence $G_{v}$ lasts at most $d s^{2}_{d}$ rounds.
\end{obse}

We now describe a strategy $S$ for Maker to obtain a candidate scheme. 
We choose $S$ in such a way that the following invariant $I$ is maintained.
\paragraph{Invariant $I$}
Suppose that $t$ rounds have been played and let $H^{\star}$ denote the corresponding position.
Let $v \in V(G)$ be a vertex which became ready in round $t$ (i.e. $v$ is ready after round $t$ but was not ready after round $t-1$). 
Then at least $s_{d}$ vertices in $S_{v}$ are untouched in $H^{\star}$. 

The invariant $I$ clearly holds for $t=0$. Indeed, every $S_{v}$ has cardinality at least $s_{d}$ and at the beginning of the game every vertex is untouched.

Let $t > 0$ and suppose that $t-1$ rounds have been played. By induction the invariant holds after round $t-1$. 
For all vertices $v$ which became ready in round $t-1$ Maker fixes a subset $S \subseteq S_{v}$ of $s_{d}$ untouched vertices and sets $B_{v} := S$. 
If all vertices of $V(G)$ are completed after round $t-1$ then by definition $(B_{v_{1}}, B_{v_{2}}, \ldots, B_{v_{n}})$ form a candidate scheme, which by Lemma \ref{lemm:CandidateSchemeImpliesGCopy} guarantees that Maker's graph contains a copy of $G$. In this case we are done. It remains to consider the case where not all vertices are completed.

Suppose that in round $t$ Breaker claims an edge $(x^{\star}, y^{\star})$ with $x^{\star} \in S_{u}$, $y^{\star} \in S_{v}$ and $l(u) < l(v)$. We distinguish three cases.
\begin{itemize}
\item[Case 1] $v$ is ready but not completed. Maker responds in the game $G_{v}$. 
\item[Case 2] $v$ is not ready. Maker selects a $w \in P(v)$ such that $w$ is ready but not completed and acts as if Breaker claimed an edge in the board of $G_{w}$ (recall that $P(v)$ denotes the set of vertices $u$ for which there is a directed path from $v$ to $u$ in $D$). Note that such a $w$ always exist. Indeed, suppose otherwise and let $w'$ be a non-ready vertex with minimum level among all vertices in $P(v)$. By construction $w'$ has out-degree at least one (otherwise $w'$ would be ready since the beginning of the game), and by assumption there is at least one out-neighbor $w''$ of $w'$ which is not completed. But then $w''$ is non-ready with $l(w'') < l(w')$, contradicting the choice of $w'$.
\item[Case 3] $v$ is completed. Maker selects a not yet completed vertex $w \in V(G)$ and acts as if Breaker claimed an edge in the board of $G_{w}$.
\end{itemize}
Note that since the boards $G_{v}$ are pairwise distinct Maker can treat each game $G_{v}$ separately.

We now show that the invariant $I$ is fulfilled after round $t$. Suppose that $v$ became ready after round $t$. We observe that for every vertex $x^{\star} \in S_{v}$ which was touched by Breaker in one of the first $t$ rounds there is (at least) one vertex $y = f(x^{\star}) \in V(G)$ such that Maker acted as if Breaker claimed an edge in the board of $G_{y}$. Note that by construction, $y \in P(v)$. 

By Observation \ref{obse:gamelengthsinglevertex}, for every $u \in P(v)$, $G_{u}$ lasted at most $d s^{2}_{d}$ rounds; thus there are at most $d s^{2}_{d}$ vertices $x^{\star} \in S_{v}$ with $f(x^{\star}) = u$. Hence there are at most $d s^{2}_{d} |P(v)|$ vertices 
$x^{\star} \in S_{v}$ where $f(x^{\star}) \in P(v)$. Therefore at most $d s^{2}_{d} |P(v)|$ vertices of $S_{v}$ have been touched during the first rounds, which yields that at least $s_{d}$ vertices remain untouched, as claimed.


Thus, if Maker follows $S$ then by Corollary \ref{coro:winningstrategyforgameononevert} he has a strategy to achieve 
that $(B_{v_{1}}, B_{v_{2}}, \ldots, B_{v_{n}})$ form a candidate scheme. Due to Lemma \ref{lemm:CandidateSchemeImpliesGCopy} this guarantees that Maker's graph contains a copy of $G$. This concludes the proof of Theorem \ref{theo:Makerwinsonsmallboard}.


\end{document}